\documentclass[a4paper,11pt]{article}

\usepackage{amsmath,amssymb,amsfonts,graphicx}
\usepackage{hyperref} 
\usepackage[latin1]{inputenc}
\usepackage{psfrag}
\usepackage{graphicx}

\newtheorem{thm}{Theorem}[section]
\newtheorem{prop}[thm]{Proposition}
\newtheorem{lem}[thm]{Lemma}

\newtheorem{rem}[thm]{Remark}

\def\be#1 {\begin{equation} \label{#1}}
\newcommand{\ee}{\end{equation}}
\def\dem {\noindent {\bf Proof : }}

\newcommand{\mb}{\medskip\noindent}
\newcommand{\gb}{\bigskip\noindent}
\newcommand{\R}{\mathbb R}

\newcommand{\C}{\mathcal C}

\newcommand{\PPP}{\mathrm P}
\def \NN {\mathrm{N}}
\newcommand{\wsto}{\stackrel{\star}{\relbar\joinrel\rightharpoonup}}
\def \virg {\, , \,\,}
\def \dsp {\displaystyle}
\def \vsp {\vspace{6pt}}

\def \I { \{1,..,p\} }

\def\sqw{\hbox{\rlap{\leavevmode\raise.3ex\hbox{$\sqcap$}}$%
\sqcup$}}
\def\findem{\ifmmode\sqw\else{\ifhmode\unskip\fi\nobreak\hfil
\penalty50\hskip1em\null\nobreak\hfil\sqw
\parfillskip=0pt\finalhyphendemerits=0\endgraf}\fi}

\title{Existence results for non-smooth second order differential inclusions, convergence result for a numerical scheme and application to the modelling of inelastic collisions}
\date{March 9, 2010}
\author{ Fr\'ed\'eric Bernicot\\ CNRS - Universit\'e Lille 1 \\ Laboratoire Paul Painlev\'e \\ 59655 Villeneuve d'Ascq Cedex, France
\\frederic.bernicot@math.univ-lille1.fr \and Aline Lefebvre-Lepot \\ CNRS - Ecole Polytechnique \\ CMAP \\
91128 Palaiseau Cedex, France\\ aline.lefebvre@polytechnique.edu }

\begin {document}
\maketitle

\begin{abstract} We are interested in existence results for second order differential inclusions, involving finite number of unilateral constraints in an abstract framework. These constraints are described by a set-valued operator, more precisely a proximal normal cone to a time-dependent set. In order to prove these existence results, we study an extension of the numerical scheme introduced in \cite{Maury} and prove a convergence result for this scheme.
\end{abstract}

\mb {\bf Key-words:} Second order differential inclusions ; Proximal normal cone ; Inelastic collisions ; Numerical scheme.

\mb {\bf MSC:} 34A60 ; 34A12 ; 65L20.


\section{Introduction}

We consider second order differential inclusions, involving proximal normal cones. These ones were firstly treated by M. Schatzman \cite{Sc} in the framework of elastic impacts and later by J.J. Moreau \cite{Moreau1, Moreau2} to model inelastic impacts for a mechanical system in order to describe contact dynamics. The impact law describing the dynamics leads to a non-increasing kinetic energy at impacts. These second order problems appear in several models of mechanical systems with a finite number of degrees of freedom and dealing with frictionless and inelastic contacts.

\mb Let us specify this class of problems. Let $I$ be a bounded time-interval, $f:I\times \R^d \rightarrow \R^d$ be a map and $C:I\rightrightarrows \R^d$ be a multi-valued map. The main question concerns the existence for solutions to the following second order differential inclusion:
\begin{equation} \label{Pb:}
\left\{
\begin{array}{l}
\dsp \forall t\in I, \quad q(t)\in C(t) \vsp \\
\dsp \frac{d^2 q}{dt^2} + \NN(C(\cdot), q(\cdot)) \ni f(\cdot, q(\cdot))  \vsp \\
\dsp \forall t\in I, \quad \dot q(t^+) = \PPP_{\C_{t,q(t)}}\dot q(t^-)\vsp \\
\dsp q(0)=q_0 \in \textrm{int}[C(0)] \vsp \\
\dsp \dot q (0)=u_0.
\end{array}
\right.
\end{equation}
We denote by $\textrm{int}[C(0)]$ the interior of the set $C(0)$, by $\NN$ the proximal normal cone and for $q\in C(t)$, by $\C_{t,q}$ the set of feasible velocities:
\be{defCq}  \C_{t,q}:=\left\lbrace u,\ q+\epsilon u \in C(t+\epsilon) \textrm{ for small enough $\epsilon>0$ }  \right\rbrace.\ee
We refer the reader to \cite{Clarke} and \cite{Clarke2} for details concerning different normal cones (``limiting cone'', ``Clarke cone'', ...). Here we will deal with ``uniform prox-regular sets'' $C(t)$ so, according to \cite{PRT}, all these cones coincide.

\begin{rem} \label{rem1}
We are looking for solutions $q$ such that $\dot q$ has a bounded variation, in order that the second order differential equation in~(\ref{Pb:}) should be thought  in the distributional sense. More precisely, we will solve it for time-measure $\dot{q}\in BV(I)$ and it should be written with time-measures 
$$  d \dot{q}  + \NN(C(\cdot), q(\cdot)) dt \ni f(\cdot, q(\cdot))dt.$$
In all this work, the second order differential inclusion will be written in the distributional sense for easiness. However, we emphasize that we consider time-measures.
\end{rem}

\mb This differential inclusion can be thought as follows: the point $q(t)$, submitted to the external force $f(t,q(t))$, has to live in the set $C(t)$ and so to follow its time-evolution. The unilateral constraint ``$q(t)\in C(t)$'' may lead to some discontinuities for the velocity $\dot q$. For example, frictionless impacts can be modelled by a second order differential inclusion involving the proximal normal cone (see \cite{Moreau1, Moreau2}). This differential inclusion does not uniquely define the evolution of the velocity during an impact. To complete the description, we impose the impact law 
$$\dot q(t^+) = \PPP_{\C_{t,q(t)}}\dot q(t^-),$$
introduced by J.J. Moreau in \cite{Moreau1} and justified by L. Paoli and M. Schatzman in \cite{Paoli-Scha2, Paoli-Scha4} (using a penalty method) for inelastic impacts.

\mb
The set $C(t)$ corresponds to a set of ``admissible configurations'' for $q$. In physical problems, it is generally described by several constaints $(g_i)_{i}$ as follows~:
\be{C:inter} C(t) := \bigcap_{i=1}^p \left\{q,\  g_i(t,q)\geq 0 \right\}. \ee

\gb
 The existence of a solution for such second-order problems is still open in a general framework. The first positive results were obtained by M.P.D. Monteiro Marques \cite{Marques} and L. Paoli and M. Schatzman \cite{Paoli-Scha3} in the case of a smooth time-independent admissible set (which locally corresponds to the single constraint case $p=1$ in (\ref{C:inter})). The proof relies on a numerical method involving a time-discretization of (\ref{Pb:}) in order to compute approximate solutions and is based on the study of its convergence . The multi-constraint case with analytical data was then treated by P. Ballard with a different method in \cite{Ballard}, where a positive result of uniqueness for such problems was obtained. Then in \cite{Paoli, Paoli2, Paoli3}, an existence result is proved in the case of a non-smooth time-independent convex admissible set (given by multiple constraints). There, the active constraints  are supposed to be linearly independent in the following sense: for each configuration $q\in \partial C$, the gradients $(\nabla g_i(q))_{i\in I}$ associated to active constraints $I:=\{i,\ g_i(q)=0\}$ are supposed to be linearly independent.

\gb In the case of non-convex admissible sets, some results were obtained for a single constraint $p=1$ (for example in \cite{DMP} or in \cite{MP} and \cite{Sc3} for the first result concerning time-dependent constraints). Recently in \cite{Maury}, B. Maury has proposed a numerical scheme for time-independent multiple and convex constraints $g_i$. The admissible set $C$ is not supposed to be convex, however at each time step, the numerical scheme uses a local convex approximation of $C$. This improvement is interesting as it permits to define an implementable scheme, since the projection onto a convex set can be performed with efficient algorithms. \\
 A first result of convergence for this scheme was proved in \cite{Maury} for a single constraint and applications to the numerical simulation of sytems of particles submitted to inelastic collisions are studied in \cite{Aline} by the second author.

\gb We emphasize that in the previously mentioned works, the different numerical schemes (permitting to discretize (\ref{Pb:})) are written (or can be written) in a multi-constraint case. The main difficulties consist in proving in a one hand the existence of solutions for (\ref{Pb:}) and in the other hand a convergence result for the associated numerical schemes for such multi-constrained problems. Concerning the uniqueness, we know from \cite{Sc} and \cite{Ballard} that even with smooth data the uniqueness does not hold. The only positive results are proved in \cite{Sc2} for one-dimensional impact problems and in \cite{Ballard} in the context of analytic data. This critical question of uniqueness is not studied here. 

\gb {\bf The framework} 

\mb In this work, we are interested in extending the previous work \cite{Maury}, in order to prove the existence of solutions and to get a convergence result of the scheme in the case of multiple time-dependent constraints. Moreover we give applications in modelling inelastic collisions.

\mb First of all, let us precise some notations. We write $W^{1,\infty}(I,\R^d)$ (resp. $W^{1,1}(I,\R^d)$) for the Sobolev space of functions in $L^\infty(I,\R^d)$ (resp. $L^1(I,\R^d)$) whose derivative is also in $L^\infty(I,\R^d)$ (resp. $L^1(I,\R^d)$). $BV(I,\R^d)$ is the space of functions in $L^\infty (I,\R^d)$ with bounded variations on $I$. We define the dual space ${\mathcal M} (I) := ({\mathcal C}_c(I))'$ where ${\mathcal C}_c(I)$ is the space of continuous functions with compact support (corresponding to the set of Radon measure due to Riesz Theorem). We set ${\mathcal M}_+ (I)$ for the subset of positive measures.

\mb We consider second-order differential inclusions involving a set-valued map $Q:[0,T]\rightrightarrows \R^d$ satisfying that for every $t\in [0,T]$, $Q(t) $ is the intersection of complements of smooth convex sets.
Let us first specify the set-valued map $Q$. -- This general framework has already been described by J. Venel in \cite{Juliette} for first order differential inclusions (fitting into the so-called sweeping process theory) and in \cite{Juliette2} for a stochastic perturbation of such problems. --\\ 
For $i \in \I$, let $g_i: [0,T] \times \R^d \rightarrow \R$ be a convex function with respect to the second variable.
For every $t \in [0,T]$, we introduce the sets $Q_i(t)$ defined by:
\begin{equation}
 Q_i(t):=\left\{ q \in \R^d \virg g_i(t,q) \geq 0 \right\},
\label{def:Qi}
\end{equation}
and the feasible set $Q(t) $ (supposed to be nonempty)
\begin{equation}
 Q(t):=\bigcap_{i=1}^{p} Q_i(t).
\label{def:Q}
\end{equation}
We denote by $I=[0,T]$ the time interval. The considered problem is the following one: we are looking for a solution
$q\in W^{1,\infty}(I,\R^d)\virg \dot q\in BV(I,\R^d)$ such that
\begin{equation} \label{Pb:cont}
\left\{
\begin{array}{l}
\dsp \forall t\in I, \quad q(t)\in Q(t) \vsp \\
\dsp \frac{d^2 q}{dt^2} + \NN(Q(\cdot), q(\cdot)) \ni f(\cdot, q(\cdot)) \vsp \\
\dsp \forall t\in I, \quad \dot q(t^+) = \PPP_{\C_{t,q(t)}}\dot q(t^-)\vsp \\
\dsp q(0)=q_0 \in \textrm{int}[Q(0)] \vsp \\
\dsp \dot q (0)=u_0,
\end{array}
\right.
\end{equation}

\noindent where $\NN(Q(t), q(t))$ is the proximal normal cone of $Q(t)$ at $q(t)$ and $\C_{t,q}$ is the set of {\it
admissible velocities} : \be{def:Cq}  \C_{t,q}:=\left\lbrace u,\ \partial_t g_i(t,q) + \langle \nabla_q \, g_i
(t,q),u\rangle\geq 0\ \textrm{ if } \  g_i(t,q)=0 \right\rbrace,\ee
which corresponds to (\ref{defCq}) in our framework.

\mb We have to make assumptions on the constraints $g_i$. First we require some regularity: we suppose that there exist $c>0$ and open sets $U_i(t) \supset Q_i(t) $ for all $ t$ in $[0,T]$  verifying
\begin{equation}
 \tag{A0}
d_H(Q_i(t), \R^d \setminus U_i(t)) > c,
\label{Ui}
\end{equation}
where $d_H $ denotes the Hausdorff distance.
Moreover we assume that there exist constants $\alpha, \beta, M >0$ such that for all $ t$ in $[0,T]$, $g_i(t,\cdot)$ belongs to $ C^2(U_i(t))$ and satisfies
\begin{equation}
 \tag{A1}
\forall \, q \in U_i(t) \virg \alpha \leq |\nabla_{q}\, g_i (t,q) | \leq \beta,
\label{gradg}
\end{equation}
\begin{equation}
 \tag{A2}
\forall \, q \in U_i(t) \virg|\partial_{t} g_i (t,q) | \leq \beta,
\label{dtg}
\end{equation}
\begin{equation}
 \tag{A3}
\forall \, q \in U_i(t) \virg  |\partial_{t}\nabla_{q}\, g_i (t,q) | \leq M,
\label{dtgradg}
\end{equation}
\begin{equation}
 \tag{A4}
\forall \, q \in U_i(t) \virg  |\mathrm{D}_q^2 g_i (t,q) | \leq M,
\label{hessg}
\end{equation}
and
\begin{equation}
 \tag{A5}
\forall \, q \in U_i(t) \virg  |\partial^2_t g_i (t,q) | \leq M.
\label{dtdtg}
\end{equation}

\mb In comparison with \cite{Juliette} and \cite{Juliette2} where first order differential inclusion are studied, we require the new and natural assumption (\ref{dtdtg}), due to the fact that we consider second order differential inclusions.

\mb
Note that these assumptions can slightly be weakened. Indeed, the lower bound in (\ref{gradg}) is only required in a neighborhood of $q\in \partial Q(t)$. Moreover, we have assumed $C^2$-smoothness in (\ref{hessg}) and (\ref{dtdtg}) for the sake of simplicity, but we only need $C^{1+\epsilon}$ regularity. 

\mb Furthermore, we require a kind of independence for the active gradients. For all $t\in [0,T]$ and $q \in Q(t)$, we denote by $I(t,q)$ the active set at $q$
\begin{equation}
 I(t,q):=\left\{i \in\I \virg  g_i(t,q)= 0 \right\},
\label{def:I}
\end{equation}
corresponding to the active constraints. For every $\rho >0 $, we define the following set:
\begin{equation}
 I_\rho(t,q):=\left\{i \in\I \virg  g_i(t,q) \leq \rho \right\}.
\label{def:Irho}
\end{equation}
We assume there exist $\gamma >0 $ and $\rho >0 $ such that for all $t \in [0,T] $,
\begin{equation}
 \tag{A6}
\forall \, q \in Q(t) \virg \forall \, \lambda_{i} \geq 0, \sum_{i \in I_\rho(t,q)} \lambda_{i} | \nabla_q \, g_i(t,q)|  \leq  \gamma\left| \sum_{i \in I_\rho(t,q)}  \lambda_{i} \nabla_q \, g_i (t,q) \right|.
\label{inegtrianginverserho}
\end{equation}
We will use the following weaker assumption too:
\begin{equation}
 \tag{A6'}
\forall \, q \in Q(t) \virg \forall \, \lambda_{i} \geq 0, \sum_{i \in I(t,q)} \lambda_{i} | \nabla_q \, g_i(t,q)|  \leq  \gamma\left| \sum_{i \in I(t,q)}  \lambda_{i} \nabla_q \, g_i (t,q) \right|.
\label{inegtrianginverse}
\end{equation}

\mb 
Note that assumptions (\ref{inegtrianginverserho}) and (\ref{inegtrianginverse}) describe a kind of ``positive linear independence'' of the almost active gradients. In the time-independent case, (\ref{inegtrianginverse}) is lightly weaker than the linear independence assumption made in~\cite{Paoli, Paoli2, Paoli3}. In fact, such assumptions imply a ``uniform prox-regularity'' of the admissible set $Q(t)$, which is a weaker property than the convexity.

\mb Under these assumptions, we have a characterization of the proximal normal cone $\NN(Q(t), \cdot)$.

\begin{prop}[Prop 2.8 of \cite{Juliette}] \label{prop:cone} In this framework, we know that for every $t\in I$, and every $q\in\partial Q(t)$,
$$ \NN(Q(t), q) := \left\{ -\sum_{i\in I(t,q)} \lambda_i \nabla_q \, g_i (t,q), \ \lambda_i\geq 0 \right\}.$$
\end{prop}
So our Problem (\ref{Pb:cont}) can be written as follows: we are looking for solutions $q\in W^{1,\infty}(I,\R^d)\virg \dot q\in BV(I,\R^d)$ and time-measures $\lambda_i \in {\cal
  M}_+(I)$ such that
\be{Pb:cont2}
\left\{
\begin{array}{l}
\dsp \forall t\in I,\quad q(t)\in Q(t) \vsp \\
\dsp \frac{d^2 q}{dt^2} = f(\cdot, q(\cdot)) + \sum_{i=1}^p \lambda_{i} \nabla_q \, g_i (\cdot ,q(\cdot))  \vsp \\
\dsp \hbox{supp}(\lambda_{i})\subset \{t\virg g_{i}(t,q(t))=0\}\hbox{ for all } i \vsp\\
\dsp \forall t\in I,\quad \dot q(t^+) = \PPP_{C_{t,q(t)}}\dot q(t^-)\vsp \\
\dsp q(0)=q_0 \in \textrm{int}[Q(0)] \vsp \\
\dot q (0)=u_0.
\end{array}
\right.
\ee

\mb We denote by $\lambda=(\lambda_1, \cdots , \lambda_p)\in\R^{p}$ the vector of the Lagrange multipliers associated to these $p$ constraints.

\mb As usual, we obtain existence results for (\ref{Pb:cont2}) by proving the convergence of a sequence of discretized solutions. \\
To do so, we extend the algorithm proposed by B. Maury in \cite{Maury} for modelling inelastic collisions, to the case of abstract and time-dependent constraints. In \cite{Maury}, the convergence (up to a subsequence) is proved in the case of a single constraint. Here, we show that this convergence still holds in the multi-constraint case. 

\mb
Let us describe the numerical scheme.
\medskip
Let $h=T/N$ be the time step. We denote by  $q_h^n\in\R^{d}$ and $u_h^n\in\R^{d}$ the approximated solution and velocity at time $t_h^nh$ for $n\in\{0,..,N\}$.

\mb
The discretization of the continuous constraints $\C_{t_h^n,q(t_h^n)}$ proposed in \cite{Maury} corresponds to a first order approximation of the constraints in the space variable: for $t\in I$ and $q\in U(t)$, we set
\be{def:Kh} K_h(t,q):=\left\lbrace u,\ g_i(t,q) + h\langle \nabla_q \, g_i (t,q),u\rangle\geq 0 \right\rbrace.\ee
The reason why we do not expand the time variable in this discrete admissible set is that we will use a semi-implicit numerical scheme, with directly impliciting in time the set $K_h(t,q)$ (See (\ref{unp1})).

\mb The approximated solutions are built using the following scheme:
\begin{enumerate}
\item Initialization : \be{CondInit} (q_h^0,u_h^0):=(q_0,u_0) \ee \item Time iterations: $q_h^n$ and $u_h^n$ are given.
We define $\dsp f_h^n:=\frac{1}{h}\int_{t_h^n}^{t_h^{n+1}}f(s,q_h^n)ds$, \be{unp1} u_h^{n+1} := \PPP_{
K_h(t_h^{n+1},q_h^{n}) } [u_h^n+hf_h^n] \ee and \be{qnp1} q_h^{n+1}:=q_h^n+hu_h^{n+1},\ee
\end{enumerate}
where $\PPP_C$ is the Euclidean projection onto the set $C$. This algorithm is a ``prediction-correction algorithm'': the predicted velocity $u_h^n+hf_h^n$, that may not be admissible, is projected onto the approximate set of admissible velocities. \\
Since the projection $\PPP_{ K_h(t_h^{n+1},q_h^{n}) }$
consists in a constrained minimization problem, with a finite number of affine constraints, it involves Lagrange
multipliers $(\lambda_{h}^{n+1})\in\R^{p}$ corresponding to the $p$ constraints. It can be checked that we have a discrete counterpart of the momentum balance appearing in~(\ref{Pb:cont2}): 
\be{eq:PFD_disc} \frac{u_h^{n+1}-u_h^n}{h} = f_h^n + \sum \lambda_{h,i}^{n+1} \nabla_q
\, g_i (t_h^{n+1},q_h^n) \ee with $\lambda_{h,i}^{n+1}\geq 0 $ and $\lambda_{h,i}^{n+1}=0$ when $g_i(t_h^{n+1},q_h^n)+h\langle \nabla_q\, g_i(t_h^{n+1},q_h^n),u_h^{n+1}\rangle>0$.

\mb In \cite{Maury}, the scheme is shown to be stable, robust and to present a good behaviour for large time-steps. That is why, we are interested in continuing its numerical analysis in the multi-constraint case, with proposing some extensions like the time-dependence of the constraints.

\gb {\bf Results}  

\mb We recall that $I=[0,T]$ is the time interval and $h$ is the constant time step ($t_h^nh$ for $n=0\ldots N$). We
denote by $q_h$ the piecewise affine function with $q_h(t_h^n)=q_h^n$. We denote by $u_h$ the derivative of $q_h$,
piecewise constant equal to $u_h^{n+1}$ on $]t_h^n,t_h^{n+1}[$. Finally, we define $\lambda_h$, piecewise constant
equal to $\lambda_h^{n+1}$ on $]t_h^n,t_h^{n+1}[$.

\mb The convergence theorem is the following one.
\begin{thm} \label{thm}
Let ($q_h,u_h,\lambda_h)$ be the sequence of solutions constructed from the scheme~(\ref{CondInit}-\ref{eq:PFD_disc}) and suppose that $f:I \times \R^d \rightarrow \R^d $ is a measurable map satisfying:
\begin{align}
  & \exists K_L >0 \virg \forall t\in I\virg \forall q,\tilde{q} \in U(t) \virg | f(t,q) -f(t,\tilde{q})| \leq K_L |q - \tilde{q}| &
\label{flip} \vsp \\
& \exists F \in L^1(I) \virg \forall t\in I\virg \forall q\in U(t) \virg |f(t,q)|\leq F(t).
\label{lingro}
\end{align}
Then, when $h$ goes to zero, there
exist subsequences, still denoted by $(q_h)_h$, $(u_h)_h$, $(\lambda_h)_h$, and
$$(q,u,\lambda) \in W^{1,\infty}(I,\R^d)\times BV(I,\R^d) \times {\mathcal M}_+(I)^p $$
such that
$$
\begin{array}{l}
u_h\longrightarrow u \hbox{ in } L^1(I,\R^d), \vsp\\
q_h\longrightarrow q \hbox{ in } W^{1,1}(I,\R^d) \hbox{ and } L^\infty(I,\R^d) \hbox{
  with } \dot q =u, \vsp\\
\lambda_h\wsto\lambda \hbox{ in } {\mathcal M}_+(I)^p
\end{array}
$$
where $(q,u,\lambda)$ is solution to~(\ref{Pb:cont2}) and so $(q,u)$ is a solution to (\ref{Pb:cont}). \label{thm:cv}
\end{thm}

\mb We emphasize that, up to our knowledge, this result is the first one concerning such multi-constrained second order differential inclusions with on the one hand time-dependent constraints and on the other hand a non-convex and non smooth admissible set.

\begin{rem} For time-independent constraints, Assumption (\ref{inegtrianginverse}) is required but Assumption (\ref{inegtrianginverserho}) is not necessary.
\end{rem}

\mb The proof is quite long and technical. We refer the reader to \cite{Maury} for a first proof dealing with the case of one ($p=1$) time-independent constraint $g$. We will follow the same reasoning with some new arguments (appearing in \cite{Juliette}) in order to solve the difficulties raised by the multiple constraints and the time-dependence. Section \ref{sec:dem} is devoted to the outline and the main ideas of the proof. For the sake of readibility, the demonstrations of some technical Propositions are postponed to Section \ref{sec:technique}. In Section \ref{sec:model}, we describe an
application to the modelling of inelastic collisions.

\section{Convergence result} \label{sec:dem}

This section is devoted to the proof of Theorem \ref{thm}. It is divided in 7 steps and for readibility reasons we have postponed some technical proofs in the next section.\\

\gb {$\bullet$ \bf Step 1: The scheme is well-defined and produces feasible configurations}

\begin{prop}
 \label{prop:feasible_config} For a small enough parameter $h$, the scheme is well-defined. Moreover the computed configurations are feasible~:
$$
\forall h>0, \forall n\in\{0,..,N\},\quad  q_h(t_h^n)\in Q(t_h^n).
$$
\end{prop}
\dem  Let $h$ be smaller than $c/c_0$ where $c$ and $c_0$ are given in Lemma \ref{lem:Qlip} (below stated) and Assumption (\ref{Ui}) respectively. By assuming that $q_h(t_h^n)\in Q(t_h^n)$, we also deduce that $q_h(t_h^n)\in Q(t_h^{n+1}) + c_0h \overline{B}(0,1) \subset
U(t_h^{n+1})$. Then the gradient $\nabla_q \, g_i (t_h^{n+1},q_h^n)$ is well-defined and so is the set $K_h(t_h^{n+1},q_h^n)$. The step $2$ of the scheme can be performed and due to
the convexity of function $g_i(t_h^{n+1},\cdot)$,
$$ u_h^{n+1}\in K_h(t_h^{n+1},q_h^n) \Longrightarrow q_h^{n+1}\in Q(t_h^{n+1}).$$
Then we conclude by iteration.
\findem

\begin{lem} \label{lem:Qlip} The set-valued map $Q$ is Lipschitz continuous with a constant $c_0$, for the Hausdorff distance.
 \end{lem}
 
\mb We refer the reader to Proposition 2.11 of \cite{Juliette} for a detailed proof of this result.

\mb For the intermediate times $t\in ]t_h^n,t_h^{n+1}[$, the point $q_h(t)$ may not belong to $Q(t)$. However from
Proposition \ref{prop:feasible_config} and Lemma \ref{lem:Qlip}, we have the following estimate~:
\be{eq:feasible_config2} \forall h>0, \forall t\in I,\quad d(q_h(t), Q(t)) \leq \max\{ d(q_h^{n},Q(t)),d(q_h^{n+1},Q(t))\} \leq c_0 h. \ee 

\gb {$\bullet$ \bf Step 2: $u_h$ is bounded in $BV(I,\R^d)$}

\mb First, we check that the velocities are uniformly bounded (proved later in Subsection \ref{subsec:uBV}).

\begin{prop} \label{prop:uLI}
 The sequence of computed velocities $(u_h)_h$ is bounded in $L^\infty(I,\R^d)$. We set 
$$ K:= \sup_{h} \|u_h\|_{L^\infty(I)} <\infty.$$ 
\end{prop}

\begin{prop}
\label{prop:uBV} The sequence of computed velocities $(u_h)_h$ is bounded in $BV(I,\R^d)$.
\end{prop}

\dem Since $u_h^0=u_0$ and using Proposition~\ref{prop:uLI}, it suffices to show that $(u_h)_h$ has bounded variations on $I$. This has been proved for a single constraint in \cite{Maury}. Unfortunately, this proof cannot be extended to the multi-constraint case. To obtain an estimate on the total variation, we use a similar technique to the one proposed in \cite{Monteiro} and \cite{DMP}. These ideas rest on the following property: all the cones $K_h(t_h^{n+1},q_h^n)$ contain a ball of fixed radius with a bounded center, which describes the fact that the solid angles of the cones $\NN(Q(t),q_h(t))$ are not too small. This property is proved by using a ``good direction'' (see Lemma \ref{lem:gooddir}) which permits to increase all the almost active constraints. The details of the proof are postponed to Subsection~\ref{subsec:uBV}. \findem

\mb{$\bullet$ \bf Step 3: Extraction of convergent subsequences}

\mb
Proposition~\ref{prop:uBV} directly implies the following convergence result:
\begin{prop}
\label{prop:q_u_cv}
 There exist $q$ in $W^{1,\infty}(I,\R^d)$ and $u$ in $BV(I,\R^d)$ such that, up to a subsequence,
$$
\begin{array}{l}
u_h\xrightarrow[h\to 0]{} u \hbox{ in } L^1(I,\R^d), \vsp\\
q_h \xrightarrow[h\to 0]{} q \hbox{ in } W^{1,1}(I,\R^d) \hbox{ and } L^\infty(I,\R^d) \hbox{
  with } \dot q =u.
\end{array}
$$
Furthermore (\ref{eq:feasible_config2}) yields
$$
\forall t\in I, \quad q(t)\in Q(t).
$$
\end{prop}
In addition, we show that the sequence of Lagrange multipliers converges:
\begin{prop}
\label{prop:lambda_cv} There exists $\lambda$ in ${\mathcal M}_+(I)^p$ such that, up to a subsequence,
$$\lambda_h\wsto\lambda \hbox{ in } {\mathcal M}_+(I)^p.$$
\end{prop}
\dem
From~(\ref{eq:PFD_disc}), we have
$$
\sum_i h\lambda_{h,i}^{n+1} \nabla_q \, g_i (t_h^{n+1},q_h^n) = u_h^{n+1}-u_h^n-hf_h^n$$
 with $\lambda_{h,i}^{n+1}\geq 0$ and $\lambda_{h,i}^{n+1}=0$ when 
 \be{eq:multiplier} g_i(t_h^{n+1},q_h^n)+h\langle \nabla_q\, g_i(t_h^{n+1},q_h^n),u_h^{n+1} \rangle>0. \ee
Remark that for a small enough parameter $h$, $g_i(t_h^{n+1},q_h^n)+h\langle \nabla_q\, g_i(t_h^{n+1},q_h^n),u_h^{n}+hf_h^n\rangle\leq0$ implies that $i\in I_\rho(t_h^{n+1},q_h^n)$.
Consequently, the reverse triangle inequality (Assumption (\ref{inegtrianginverse})) with (\ref{gradg}) and the Lipschitz regularity (Assumption (\ref{dtgradg})) imply for a small enough parameter $h$
$$ h\sum_{i} \lambda_{h,i}^{n+1}  \lesssim |u_h^{n+1}-u_h^n-hf_h^n|.$$
Therefore, we obtain for this small enough parameter $h$
$$
\|\lambda_h\|_{L^1(I)}\lesssim \textrm{Var}_{I}(u_h) + \|F\|_{L^1(I)}.
$$
Hence from Proposition~\ref{prop:uBV} together with hypothesis~(\ref{lingro}), $(\lambda_h)_h$ is bounded in $L^1(I)$, which concludes the proof. \findem

\mb{$\bullet$ \bf Step 4: Momentum balance}

\mb
As in Step 6 of Theorem 1 in \cite{Maury}, using Proposition~\ref{prop:q_u_cv} together with Proposition~\ref{prop:lambda_cv}, we pass to the limit in the discrete momentum balance~(\ref{eq:PFD_disc}) to obtain
\begin{prop}
The momentum balance is verified by the limits $u$ and $\lambda$: in the sense of time-measure
$$\dot u =  f(\cdot,q(\cdot)) + \sum_{i}\lambda_{i} \nabla_q \, g_i (\cdot,q(\cdot)).$$
\end{prop}
Note that this equation has to be thought in term of time-measure (see Remark \ref{rem1}):
$$ du =  f(\cdot)d t + \sum_{i} \nabla_q \, g_i (\cdot,q(\cdot)) \lambda_{i},$$
where we denote by $du$ the differential measure of the $BV$-function $u$.

\mb{$\bullet$ \bf Step 5: Support of the measures $\lambda_{i}$}

\mb
From the uniform convergence of $q_h$ and the Lipschitz regularity of $g_i$ (Assumptions (\ref{gradg}) and (\ref{dtg})), it can be checked, as in Step 7 of Theorem 1 in \cite{Maury}, that
\begin{prop}
$$ \forall i\virg \hbox{supp}(\lambda_{i})\subset \{t\virg g_i(t,q(t))=0\}.$$
\end{prop}

\mb Indeed (\ref{eq:multiplier}) describes a similar property for the discretized multipliers. The uniform convergence allows us to go to the limit in (\ref{eq:multiplier}) and to prove the previous proposition. \\
This property describes the fact that the measure $\lambda_i$ has a contribution only when the associated constraint $g_i$ is saturated. 

\mb{$\bullet$ \bf Step 6: Initial condition}

\mb
As in Step 8 of Theorem 1 in \cite{Maury}, using again the uniform convergence of $q_h$ it can be shown that $$q(0)=q_0 \qquad \textrm{and} \qquad u(0)=u_0.$$
We emphasize that to prove this point, we use the property $q_0\in \textrm{Int}[Q(0)]$. From this, it can be checked that for $t_h^n <s$ with $s$ a small enough parameter the desired velocity $u_h^n+hf_h^n$ still remains admissible and so we don't need to project. 
That allows us to deal with any initial velocity $u_0\in \R^d$. If $q_0 \in \partial Q(0)$, this property still holds if we assume that the initial velocity is admissible: $u_0\in {\mathcal C}_{0,q_0}$. Else, we would get
$$ u^+(0)=\PPP_{\C_{0,q_0}}(u_0)$$
according to the next proposition.

\mb{$\bullet$ \bf Step 7: Collision law}

\mb Finally, Theorem~\ref{thm:cv} will follow, provided that we check the collision law for the limits $u$ and $q$,
which is given by the following proposition.
\begin{prop}
 \label{prop:CollLaw}
$$ \forall t_0\in I\virg u^+(t_0) = \PPP_{\C_{t_0,q(t_0)}}(u^-(t_0)).$$
\end{prop}
\dem The idea is to let $h$ go to zero in the discrete collision law,
$$
u_h^{n+1}=\PPP_{K_h(t_h^{n+1},q_h^{n})} [u_h^n+hf_h^n].
$$
The main difficulty comes from the fact that the mapping $q\rightarrow K_h(t,q)$ is not Lipschitzean. The details of
the proof are postponed to Subsection~\ref{subsec:CollLaw}. \findem

\section{Auxiliary results} \label{sec:technique}

\mb Before proving the technical Propositions \ref{prop:uLI}, \ref{prop:uBV} and \ref{prop:CollLaw}, we recall the following main Lemma. This technical result is very important and all our proofs rest on this idea. It says that, for each time $t_h^n$, one can find a ``good direction'' increasing all the constraints which are almost active, the corresponding increase being independent of $n$.

\begin{lem} \label{lem:gooddir} There exist constants $\delta,\kappa,\theta$ and $\tau>0$ such that for all $t\in I$, for all $q\in \partial Q(t)$ there exists a unit vector $v:=v(t,q)$ satisfying~:
\begin{itemize}
  \item for all $s\in[t-\tau,t+\tau]$, $y\in Q(s)\cap B(q,\theta)$ and  $i\in I_{\kappa \rho}(s,y)$,
$$\left<\nabla_q\, g_i(s,y),v \right>\geq\delta,$$
where $\rho$ is the constant defined in~(\ref{inegtrianginverserho}). 
\end{itemize}
\end{lem}

\mb This lemma is a consequence of the reverse triangle inequality (Assumption (\ref{inegtrianginverserho})). We do not give the proof here and refer the reader to Lemma 2.10 in \cite{Juliette} for a detailed proof with $\tau=0$ and to Proposition 4.2 in \cite{Juliette2} for a complete proof with some $\tau>0$. Indeed there is in the previously mentioned papers, a detailed construction of such ``good directions''.

\subsection{BV estimate for $u_h$ (Propositions~\ref{prop:uLI} and \ref{prop:uBV})}
\label{subsec:uBV}

We first prove a uniform bound of the computed velocities $u_h$ in $L^\infty(I)$.

\mb
{\bf Proof of Proposition \ref{prop:uLI}} \\
 {\bf $1-)$} For any $t\in I$ and $q\in Q(t)$, construction of a specific point $w\in \C_{t,q}$.\\
From Lemma \ref{lem:gooddir} (stated at Section \ref{sec:technique}),  
 it exists a ``good direction'': a unit vector $v$ satisfying, 
\be{eq:gooddirtnp1} 
\forall i\in I_{l}(t,q),\quad
\left<\nabla_q\, g_i(t,q),v \right>\geq\delta,
\ee
for some numerical constants $\delta,l>0$ non-depending on $t$ and $q$. For $k>0$ a large enough real, we claim that $kv$ belongs to $\C_{t,q}$.
Indeed for $i\in I_l(t,q)$ ($\supset I(t,q)$), with $k\geq (\beta+\delta)/\delta$, it comes
\be{eq:error} \partial_t g_i(t,q) + \langle \nabla_q \, g_i
(t,q),kv\rangle\geq \partial_t g_i(t,q) + k\delta\geq \delta>0 \ee
thanks to Assumption (\ref{dtg}). \\
Choosing $k=(\beta+\delta)/\delta$, we have built a point $w=kv$ belonging to $\C_{t,q}$ with $|w|=(\beta+\delta)/\delta$.

\mb {\bf $2-)$} This vector $w$ belongs to $K_h(s+h,\tilde{q})$ for $(s,\tilde{q})$ close to $(t,q)$ and $h$ small enough.\\
Let fix a point $(t,q)$ (for example the initial condition $(t,q)=(0,q_0)$). From the previous point, we know that there exists a bounded admissible velocity $w\in \C_{t,q}$, which satisfies the stronger property (\ref{eq:error}). \\
More precisely if $s\in I$ and $\tilde{q}\in Q(s)$ satisfy
\be{eq:error2} |s-t|+|\tilde{q}-q| \leq \min\{ l/(2\beta),\epsilon\}:=\nu \ee with $\epsilon$ a small parameter verifying $2\epsilon M(1+(\beta +\delta)/\delta)\leq \delta$  then $w\in K_h(s+h, \tilde q)$. Indeed, from (\ref{eq:error2}) we have
$$ I_{l/2}(s,\tilde{q}) \subset I_l(t,q).$$
This, together with (\ref{eq:error}) and (\ref{eq:error2}) gives, for all $i\in
I_{l/2}(s,\tilde{q})$
 $$ \partial_t g_i(s,\tilde{q}) + \langle \nabla_q\, g_i(s,\tilde{q}),w \rangle \geq \delta-\epsilon M(1+|w|)\geq \delta/2 .$$
Moreover, for the other indices $i\notin I_{l/2}(s,\tilde{q})$ we have
$$ g_i(s,\tilde{q}) + h\left[ \partial_t g_i(s,\tilde{q}) + \langle \nabla_q\, g_i(s,\tilde{q}),w \rangle\right] \geq \frac{l}{2} -h\beta(1+|w|).$$
Consequently, for $h$ small enough ($h\leq h_0:=l/(2\beta(1+(\beta+\delta)/\delta)+\delta/2)$), we obtain
$$  \forall i, \qquad g_i(s,\tilde{q}) + h\left[ \partial_t g_i(s,\tilde{q}) + \langle \nabla_q\, g_i(s,\tilde{q}),w \rangle\right] \geq h\frac{\delta}{2},$$
which, by a first order expansion in time gives:
$$  \forall i, \qquad g_i(s+h,\tilde{q}) + h\langle \nabla_q\, g_i(s+h,\tilde{q}),w \rangle \geq h\frac{\delta}{2} +O_{h\to 0}(h^2).$$
We deduce that there exists $h_1\leq h_0$ such that, for $h\leq h_1$,
$$  \forall i, \qquad g_i(s+h,\tilde{q}) + h\langle \nabla_q\, g_i(s+h,\tilde{q}),w \rangle \geq 0,$$
and consequently $w\in K_h(s+h,\tilde{q})$.

\mb {\bf $3-)$} Estimate on the velocities for small time intervals.\\
Let us fix $h\leq h_1$ (given in the previous point) and a small time interval $[t_-,t_+]\subset I$ of length 
$$ |t_+-t_-|\leq \frac{\nu}{2\left(|u_h^{n_0}|+2\frac{\beta + \delta}{\delta} + \int_0^T F(t) dt\right)}$$
where $n_0$ is the smallest integer $n$ such that $t_h^n\geq t_-$. We suppose that $t_h^{n_0}\in [t_-,t_+]$.
We are looking for a bound on the velocity on this time interval.
From the first two points, setting $(t,q)=(t_h^{n_0},q_h^{n_0})$, we have an admissible velocity $w\in \C_{t,q}$ such that for all $s\in I$ and $\tilde{q}\in Q(s)$ we have
$$ w\in K_h(s+h,\tilde{q})$$
as soon as $|s-t|+|\tilde{q}-q|\leq \nu$.
Since for $s=t_h^n\in [t_-,t_+]$ , $|s-t|\leq \nu/2$, we deduce that for all integer $n$ such that $t_h^n\in [t_-,t_+]$, if
\be{eq:error3} |q_h^n-q_h^{n_0}|\leq \nu/2 \ee then $$ w\in K_h(t_h^{n+1},q_h^n).$$
Considering such an integer $n$ satisfying~(\ref{eq:error3}), since $u_h^{n+1}$ is the Euclidean projection of $u_h^{n}+hf_h^{n}$ on the convex set $K_h(t_h^{n+1},q_h^{n})$ (containing the point $w$), we deduce that
$$ |u_h^{n+1} - w|\leq | u_h^{n}+hf_h^n - w|,$$
which implies
$$ |u_h^{n+1} - w|\leq |u_h^{n}-w| + \int_{t_h^{n}}^{t_h^{n+1}} F(t) dt.$$
We set $m$ the smallest integer (bigger than $n_0$) such that $m+1$ does not satisfy (\ref{eq:error3}) or $t_h^{m+1}\notin [t_-,t_+]$. By summing these inequalities from $n_0$ to $n=p-1$ with $n_0\leq p\leq m$, we get
$$\forall p\in[n_0,m],\qquad |u_h^{p}-w|\leq |u_h^{n_0}-w| + \int_0^T F(t) dt.$$
Finally, it comes
\be{eq:error4} \sup_{n_0\leq p\leq m} |u_h^{p}| \leq |u_h^{n_0}|+2\frac{\beta + \delta}{\delta} + \int_0^T F(t) dt.\ee
By integrating in time, we deduce
$$ |q_h^{m+1}-q_h^{n_0}| \leq \left(|u_h^{n_0}|+2\frac{\beta + \delta}{\delta} + \int_0^T F(t) dt\right)|t_+-t_-| \leq \nu/2$$
by the assumption on the length of the time interval.
As a consequence, we get that $n=m+1$ satisfies (\ref{eq:error3}) which by definition of $m$, yields $t_h^{m}\leq t_+ < t_h^{m+1}$. Hence, from~(\ref{eq:error4}), we have
$$ \sup_{t_-\leq t_h^n\leq t_+} |u_h^{n}| \leq |u_h^{n_0}|+2\frac{\beta + \delta}{\delta} + \int_0^T F(t) dt.$$

\mb {\bf $4-)$} End of the proof.\\
The parameter $h<h_1$ being fixed, we are now looking for a bound on $u_h$ on the whole time interval $I=[0,T]$. Let us start with $t_-=t(0):=0$. From the previous point we know that with 
$$ t_+=t(1):= \min\left\{\frac{\nu}{2\left(|u_0|+2\frac{\beta + \delta}{\delta} + \int_0^T F(t) dt\right)},T \right\}$$
we have
$$ \sup_{0\leq t_h^n\leq t(1)} |u_h(t_h^n)| \leq |u_0|+2\frac{\beta + \delta}{\delta} + \int_0^T F(t) dt.$$
Then, let us suppose that there exists $n_1$ such that $t(0)<t^{n_1}_h\leq t(1)<t^{n_1+1}_h$. We have $0\leq\delta_1:=t(1)-t^{n_1}_h<h$. In that case, we set $t_-=t_h^{n_1}$ and  
\begin{align*} 
t_+=t(2)&:= \min\left\{t^{n_1}+\frac{\nu}{2\left(|u_0|+4\frac{\beta + \delta}{\delta} + 2\int_0^T F(t) dt\right)},T\right\}\\
&=\min\left\{t(1)-\delta_1+\frac{\nu}{2\left(|u_0|+4\frac{\beta + \delta}{\delta} + 2\int_0^T F(t) dt\right)},T\right\}.
\end{align*}
From the previous point, we deduce that
$$ \sup_{t(1)\leq t_h^n\leq t(2)} |u_h(t_h^n)| \leq \sup_{t^{n_1}\leq t_h^n\leq t(2)} |u_h(t_h^n)| \leq |u_0|+4\frac{\beta + \delta}{\delta} + 2\int_0^T F(t) dt$$
and so
$$ \sup_{0 \leq t_h^n\leq t(2)} |u_h(t_h^n)| \leq |u_0|+4\frac{\beta + \delta}{\delta} + 2\int_0^T F(t) dt.$$
By iterating this reasoning, for any integer $k\geq 1$ we set
\begin{align*} 
t(k)& := \min\left\{t(k-1)-\delta_{k-1}+ \frac{\nu}{2\left(|u_0|+2k\frac{\beta + \delta}{\delta} + 2k\int_0^T F(t) dt\right)},T\right\} \\
 & = \min\left\{-\sum_{i=1}^{k-1}\delta_i+\sum_{i=1}^k \frac{\nu}{2\left(|u_0|+2i\frac{\beta + \delta}{\delta} + 2i\int_0^T F(t) dt\right)},T\right\}.
\end{align*}
where $\delta_k<h$ for all $k$. This construction of $t(k)$ can be made while there exists $n_k$ such that $t(k-2)<t^{n_{k-1}}_h\leq t(k-1)<t^{n_{k-1}+1}_h$. That is, while $t(k-1)-t(k-2)>h$. This condition will be verified as long as
$$ -\delta_{k-2}+ \frac{\nu}{2\left(|u_0|+2(k-1)\frac{\beta + \delta}{\delta} + 2(k-1)\int_0^T F(t) dt\right)}>h.$$
Therefore, using the fact that $0\leq \delta_{k-2}<h$, we see that we can construct $t(k)$ for $k<N$ verifying
$$\frac{\nu}{2\left(|u_0|+2(k-1)\frac{\beta + \delta}{\delta} + 2(k-1)\int_0^T F(t) dt\right)}>2h,$$
which is equivalent to
\be{k0}
k<k_0(h):=1+\left(\frac{\nu}{8h}-\frac{|u_0|}{2}\right) \left(\frac{\beta+\delta}{\delta} + \int_0^T F(t) dt \right)^{-1}.
\ee
Consequently, we know that the velocities can be bounded on $[0,t(k_0(h))]$ where
$$
t(k_0(h))=\min\left\{-\sum_{i=1}^{k_0(h)-1}\delta_i+\sum_{i=1}^{k_0(h)} \frac{\nu}{2\left(|u_0|+2i\frac{\beta + 1}{\delta} + 2i\int_0^T F(t) dt\right)},T\right\}.
$$
Now, using the fact that $k_0(h)$ goes to infinity when $h$ goes to zero, that the harmonic serie diverges and that (\ref{k0}) yields
$$
\left|\sum_{i=1}^{k-1}\delta_i\right|\leq h k_0(h)\leq C,
$$
we see that $t(k_0(h))$ is equal to $T$ for $h$ small enough. Therefore, there exists $h_2<h_1$ such that, for $h<h_2$, $T=t(k_0(h_2))=t(k_0(h))$. Finally, we see that, for $h<h_2$, $t(k)$ can be constructed until $k=k_0(h_2)$ and $u_h$ can be bounded as follows
\be{eq:final} \sup_{h\leq h_2} \sup_{0\leq t_h^n\leq T} |u_h(t_h^n)| \leq |u_0|+2k_0(h_2)\frac{\beta + \delta}{\delta} + 2k_0(h_2)\int_0^T F(t) dt \ee
which concludes the proof of the existence of a uniform bound in $L^\infty$ for the velocities $u_h$.
\findem

\mb To prove Proposition~\ref{prop:uBV}, it suffices now to show that the sequence $(u_h)_h$ has bounded variation.

\begin{thm} \label{thm:uBV}
 The sequence $(u_h)_h$ has a bounded variation on $I$.
\end{thm}

\dem
In order to study the variation of $u_h$ on $I$, we split $I$ into smallest intervals. We define $(s_j)_j$ for $j$ from $0$ to $P$ such that:
$$
\left|\begin{array}{l}
\dsp s_0=0\virg s_P=T,\smallskip\\
\dsp|s_{j+1}-s_j|=\frac{1}{2}\min\left\{\tau, \frac{\theta}{K}\right\}, \hbox{ for } j=0\ldots P-2,\smallskip\\
\dsp|s_{P}-s_{P-1}|\leq \frac{1}{2}\min\left\{\tau, \frac{\theta}{K}\right\},
\end{array}\right.$$
where $\tau$ and $\theta$ are given by Lemma~\ref{lem:gooddir} and $K$
is the bound on $\|u_h\|_{L^\infty(I)}$ (see Proposition
\ref{prop:uLI}). All these constants do not depend on $h$ and such a construction gives
\be{eq:boundP}
P=\left[\frac{2T}{\min\left\{\tau, \frac{\theta}{K}\right\}}\right]+1,
\ee 
which is independent of $h$.
Then, for all $h$, we define $n^j_h$ for $j$ from $0$ to $P-1$ as the first time step strictly greater than $s_j$: $$t_h^{n_h^j-1}\leq s_j < t_h^{n_h^j}, $$ and $n_h^P$ is set equal to $N$ ($t_h^N=t_h^{n_h^P}=T$). 

\mb
In the following, we suppose $h<\min\{|s_{j+1}-s_j|\}/2$. Doing so, we obtain a strictly increasing sequence of $(t_h^{n_h^j})_j$ with
\be{eq:dtnj}
|t_h^{n_h^j}-t_h^{n_h^{j-1}}|\leq \min\left\{\tau, \frac{\theta}{K}\right\}.
\ee

\mb The variation of $u_h$ on $I$ can be written as follows
$$
\textrm{Var}_I (u_h)=\sum_{n=0}^{N-1} |u_h^{n+1}-u_h^n| = \sum_{j=0}^{P-1} \textrm{Var}_j u_h
$$
where
$$
\textrm{Var}_j (u_h) := \sum_{n_h^j}^{n_h^{j+1}-1} |u_h^{n+1}-u_h^n|
$$
corresponds to the variation on $[t_h^{n_h^j},t_h^{n_h^{j+1}}[$. To study these terms, we recall that
\be{eq:unp1}
u_h^{n+1} = \PPP_{ K_h(t_h^{n+1},q_h^{n}) } [u_h^n+hf^n]
\ee
by construction and state the following lemma:

\begin{lem}\label{lem:x1_moins_x0}
There exist $\eta>0$ and uniformly bounded vectors $y^{n_h^j}$ such that, for all small enough $h$, for all $j=0\ldots P$ and $n\in [n_h^j,n_h^{j+1}[\virg $ we have
$$
x_1 = \PPP_{ K_h(t_h^{n+1},q_h^{n}) } [x_0]\quad \Longrightarrow \quad |x_1-x_0|\leq \frac{1}{2\eta}\left(|x_0-y^{n_h^j}|^2-|x_1-y^{n_h^j}|^2 \right)
$$
\end{lem}

\dem The outline of the proof is the following: first, we prove that there exist unit vectors $v^{n_h^j}$ such that \be{eq:vni} n\in
[n_h^j,n_h^{j+1}[\quad \Longrightarrow \quad \overline{B}(\frac{2\beta K}{\delta} v^{n_h^j},\eta) \subset K_h(t_h^{n+1},
q_h^{n})\  \hbox{ with } \eta:=\frac{K}{2}, \ee where $K$ is a bound on $\|u_h\|_{L^\infty(I)}$ (see Proposition \ref{prop:uLI}). Then, we conclude using similar arguments to the ones exposed in
\cite{Monteiro,DMP}.

\mb{\bf Step 1}: From Lemma~\ref{lem:gooddir} with $t=t_h^{n_h^j}$ and $q=q_h^{n_h^j}$, we have a unit ``good direction''
written $v^{n_h^j}$. Let $n$ belong to $[n_h^j,n_h^{j+1}[$. From Proposition \ref{prop:feasible_config}, we know that $q_h^{n+1}$ belongs to
$Q(t_h^{n+1})$. Moreover,~(\ref{eq:dtnj}) gives $|t_h^{n+1}-t_h^{n_h^j}|\leq \tau$ and $|q_h^{n+1}-q_h^{n_h^j}|\leq \theta$. Consequently, Lemma~\ref{lem:gooddir} gives
\be{eq::} \forall\, i\in I_{\kappa\rho}(t_h^{n+1},q_h^{n+1})\virg \left<\nabla_q\, g_i(t_h^{n+1},q_h^{n+1}),v^{n_h^j} \right>\geq\delta. \ee
We deduce that for all index $i\in \I$ and a small enough parameter $h$~:
\be{eq:gi} g_i(t_h^{n+1},q_h^{n}) + \frac{2\beta K}{\delta} h\langle \nabla_q\, g_i(t_h^{n+1},q_h^{n+1}),v^{n_h^j} \rangle \geq -h\beta K + 2h\beta K= h\beta K.\ee
Indeed, we write
$$
\begin{array}{l}
\dsp g_i(t_h^{n+1},q_h^{n}) + \frac{2\beta K}{\delta} h\langle \nabla_q\, g_i(t_h^{n+1},q_h^{n+1}),v^{n_h^j} \rangle = \smallskip \\ 
\dsp \quad\quad \left\{ g_i(t_h^{n+1},q_h^{n}) - g_i(t_h^{n+1},q_h^{n+1})\right\}+ \left\{g_i(t_h^{n+1},q_h^{n+1})+\frac{2\beta K}{\delta} h\langle \nabla_q\, g_i(t_h^{n+1},q_h^{n+1}),v^{n_h^j} \rangle\right\}.
\end{array}
$$
The first term can be estimated using~(\ref{gradg}) and the bound $K$ on $\|u_h\|_{L^\infty(I)}$. In order to estimate the second term, if $i\in I_{\kappa\rho}(t_h^{n+1},q_h^{n+1})$, we use (\ref{eq::}) together with the fact that $g_i(t_h^{n+1},q_h^{n+1})\geq 0$, which gives the required bound. In the case $i\notin I_{\kappa\rho}(t_h^{n+1},q_h^{n+1})$, we use $g_i(t_h^{n+1},q_h^{n+1})\geq \kappa \rho$ and~(\ref{gradg}), which also gives the required bound for $h$ small enough. \\
Finally,~(\ref{eq:gi}) together with assumption (\ref{gradg}) implies that for all $v\in \overline{B}(\frac{2\beta K}{\delta} v^{n_h^j},K/2)$
$$ g_i(t_h^{n+1},q_h^n) + h\langle \nabla_q\, g_i(t_h^{n+1},q_h^{n}),v \rangle \geq h\beta K -\frac{h\beta K}{2}  = \frac{h\beta K}{2}\geq 0,$$
which proves (\ref{eq:vni}).

\mb{\bf Step 2}:
Let $n$ belong to $[n_h^j,n_h^{j+1}[$. We define
$$
z^{n_h^j}:=y^{n_h^j}+\eta\frac{x_0-x_1}{|x_0-x_1|}
\quad\hbox{where}\quad
y^{n_h^j}:=\frac{2\beta K}{\delta} v^{n_h^j}.$$
(Here we suppose $x_0\neq x_1$, else the desired result is obvious.)
From the previous step we have
$$
z^{n_h^j}\in \overline{B}(\frac{2\beta K}{\delta} v^{n_h^j},\eta)\subset
K_h(t_h^{n+1},q_h^{n}).
$$
The point $x_1$ being the projection of $x_0$
onto the closed convex set $K_h(t_h^{n+1},q_h^{n})$, we get
$$
\langle x_0-x_1,z^{n_h^j}-x_1 \rangle \leq 0.
$$
From this we have
\begin{eqnarray*}
 |x_0-y^{n_h^j}|^2&=&|x_1-y^{n_h^j}|^2+|x_0-x_1|^2+2\langle z^{n_h^j}-y^{n_h^j},x_0-x_1\rangle+2\langle x_1-z^{n_h^j},x_0-x_1\rangle\\
&\geq&|x_1-y^{n_h^j}|^2+2\langle z^{n_h^j}-y^{n_h^j},x_0-x_1\rangle\\
&\geq&|x_1-y^{n_h^j}|^2+2\eta|x_0-x_1|.
\end{eqnarray*}
 This, together with the fact that the vectors $y^{n_h^j}$ are uniformly bounded by $\frac{2\beta K}{\delta}$, ends the proof of Lemma~\ref{lem:x1_moins_x0}. \findem

\mb We now come back to the proof of Theorem \ref{thm:uBV}. For $n$ in $[n_h^j,n_h^{j+1}[$, using~(\ref{eq:unp1}) and the previous lemma (with $x_0=u_h^n+hf_h^n$ and $x_1=u_h^{n+1}$), it comes
\begin{eqnarray*}
 |u_h^{n+1}-u_h^n-hf_h^n|&\leq&\frac{1}{2\eta}\left(|x_0-y^{n_h^j}|^2-|x_1-y^{n_h^j}|^2\right)\\
&\leq&\frac{1}{2\eta}\left(|u_h^n+hf_h^n-y^{n_h^j}|^2-|u_h^{n+1}-y^{n_h^j}|^2\right)\\
&\leq&\frac{1}{2\eta}\left(|u_h^n-y^{n_h^j}|^2-|u_h^{n+1}-y^{n_h^j}|^2\right) +\frac{1}{2\eta}|hf_h^n|^2+\frac{1}{\eta}|hf_h^n||u_h^n-y^{n_h^j}|\\
&\leq&\frac{1}{2\eta}\left(|u_h^n-y^{n_h^j}|^2-|u_h^{n+1}-y^{n_h^j}|^2\right)
+\frac{1}{2\eta}|hf_h^n|^2+\frac{1}{\eta}|hf_h^n|\left(K+L\right)
\end{eqnarray*}
where $L:=2\beta K/\delta$ (see Proposition \ref{prop:uLI} for the definition of $K$). By summing up these terms
for $n$ from $n_h^j$ to $n_h^{j+1}-1$ we get
\begin{eqnarray*}
\textrm{Var}_j (u_h) = \sum_{n_h^j}^{n_h^{j+1}-1} |u_h^{n+1}-u_h^n|&\leq& \frac{1}{2\eta}\left(|u_h^{n_h^{j}}-y^{n_h^j}|^2-|u_h^{n_h^{j+1}}-y^{n_h^j}|^2\right) +\sum_{n_h^j}^{n_h^{j+1}-1}\frac{1}{2\eta}|hf_h^n|^2\\
&&\quad\quad+\frac{1}{\eta}\left(K+L+\eta\right)\sum_{n_h^j}^{n_h^{j+1}-1}|hf_h^n|
\end{eqnarray*}
and finally
\begin{eqnarray*}
\textrm{Var} (u_h) =\sum_{j=0}^{P-1} \textrm{Var}_j (u_h) & \leq &
\frac{1}{2\eta}\sum_{j=0}^{P-1}\left(|u_h^{n_h^{j}}-y^{n_h^j}|^2-|u_h^{n_h^{j+1}}-y^{n_h^j}|^2\right)
+ \frac{1}{2\eta}\|F\|_{L^1(I)}^2\\
& &\quad\quad+\frac{1}{\eta}\left(K+L+\eta\right)\|F\|_{L^1(I)}\\
&\leq&\frac{1}{\eta}\left(K+L\right)^2 P+
\frac{1}{2\eta}\|F\|_{L^1(I)}^2+\frac{1}{\eta}\left(K+L+\eta\right)\|F\|_{L^1(I)}.
\end{eqnarray*}
This completes the proof of
Theorem~\ref{thm:uBV}, since $P$ does not depend on $h$ from~(\ref{eq:boundP}). \findem

\subsection{Collision law for the limits $u$ and $q$ (Proposition~\ref{prop:CollLaw})}
\label{subsec:CollLaw}

This subsection is devoted to the proof of Proposition \ref{prop:CollLaw}, recalled in the following Theorem~:

\begin{thm} \label{thm:CollLaw}
Let $t_0\in I$ be fixed. The limit function $u$ verifies:
$$ u^+(t_0) = \PPP_{\C_{t_0,q(t_0)}}(u^-(t_0)).$$
\end{thm}

\mb Note that, from Proposition \ref{prop:q_u_cv}, $u\in BV(I)$, so that the left-sided $u^-(t_0)$ and the right-sided $u^+(t_0)$ limits are well-defined. \\
The proof is quite technical  so for an easy reference, we remember the definitions of the sets $\C_{t,q}$ (given in (\ref{def:Cq})):
$$\C_{t,q}:=\left\lbrace u,\ \partial_t g_i(t,q)+\langle \nabla_q\, g_i(t,q),u\rangle\geq 0,\ \textrm{if} \  g_i(t,q)=0 \right\rbrace$$
and $K_h(t,q)$ (given in (\ref{def:Kh})):
$$K_h(t,q):=\left\lbrace u,\ g_i(t,q) + h\langle \nabla_q\, g_i(t,q),u\rangle\geq 0 \right\rbrace.$$
Moreover, we recall that
$$ K:= \sup_h \|u_h \|_{L^\infty(I)}<\infty.$$
The desired property
\be{eq:coll} u^+(t_0) = \PPP_{\C_{t_0,q(t_0)}}(u^-(t_0)) \ee
 can be seen as the limit (for $h$ going to $0$) of the ``discretized property''
\be{eq:projete0}
u_h^{n+1}=\PPP_{K_h(t_h^{n+1},q_h^{n})} [u_h^n+hf^n].
 \ee

\dem
First we claim that
\be{eq:u+}
u^+(t_0)\in C_{t_0,q(t_0)}.
\ee
To verify this property, let us consider an index $i$ such that $g_i(t_0,q(t_0))=0$. Then a first order expansion  gives~:
$$ g_i(t_0+\epsilon,q(t_0+\epsilon)) = \epsilon \left[\partial_t g_i(t_0,q(t_0))+\langle u^+(t_0),\nabla_q\, g_i(t_0,q(t_0)) \rangle\right] + o_{\epsilon\to 0}(\epsilon). $$
The feasibility of $q(t_0+\epsilon)$  (see Proposition \ref{prop:q_u_cv}) yields
$$
\partial_t g_i(t_0,q(t_0)) + \langle u^+(t_0),\nabla_q\, g_i(t_0,q(t_0)) \rangle\geq 0$$
which corresponds to (\ref{eq:u+}). \\
Let us now come back to the proof of~(\ref{eq:coll}). As we just proved $u^+(t_0)\in \C_{t_0,q(t_0)}$ and since
$\C_{t_0,q(t_0)}$ is a convex set, (\ref{eq:coll}) is equivalent to 
\be{amontrer} \forall w\in \C_{t_0,q(t_0)}, \qquad
\langle u^-(t_0) -u^+(t_0), w-u^+(t_0) \rangle \leq 0.\ee
So, in the following, let us choose $w\in \C_{t_0,q(t_0)}$. To prove (\ref{amontrer}), we construct a family of points $w_\nu$ for $\nu>0$ such that $w_\nu$ tends to $w$ when $\nu$ goes to zero and satisfies
$w_\nu\in K_h(t+h,q)$ for $h$ sufficiently small and $(t,q)$ close to $(t_0,q(t_0))$. Then, for each $\nu$, we go to the limit  on $h$, $t$ and $q$ to show that $\langle u^-(t_0)-u^+(t_0), w_\nu-u^+(t_0) \rangle \leq 0$ and finally, we make $\nu$ go to zero to conclude.

\mb{\bf Step 1}: From Lemma \ref{lem:gooddir}, there exists a neighborhood $U \subset I \times \R^d$ around $(t_0,q(t_0))$ and $v\in \R^d$ such that for all $t\in I$ and $q\in
Q(t)$
\begin{equation}
(t,q)\in U \quad\Longrightarrow\quad \forall i\in I_{\kappa\rho}(t,q) \virg \langle \nabla_q\,  g_i(t,q),v \rangle \geq \delta, \label{eq:gooddir}
 \end{equation}
with a numerical constant $\delta>0$. For $\nu>0$, we consider the point $w_\nu:=w+\nu v$ with $\nu> 0$. For all $i\in
I_{\kappa\rho}(t,q) \cap I(t_0,q(t_0))$,~(\ref{eq:gooddir}) together with $w\in\C_{t_0,q(t_0)}$ gives
\begin{align*}
 \partial_t g_i(t_0,q(t_0))+ \langle \nabla_q\, g_i(t,q),w_\nu \rangle & = \partial_t g_i(t_0,q(t_0))+\langle \nabla_q\, g_i(t,q),w \rangle + \nu \langle \nabla_q\, g_i(t,q),v \rangle \\
 & \geq \langle \nabla_q\, g_i(t,q)-\nabla_q\, g_i(t_0,q(t_0)),w \rangle + \nu\delta \\
 & \geq \nu\delta - M |w| \left[|t-t_0|+|q-q(t_0)|\right]
\end{align*}
and consequently from Assumptions (\ref{dtgradg}) and (\ref{dtdtg})
$$\partial_t g_i(t,q)+ \langle \nabla_q\, g_i(t,q),w_\nu \rangle \geq \nu\delta - (M |w|+M)
\left[|t-t_0|+|q-q(t_0)|\right].$$
 So for every $\nu>0$, if $(t,q)$ is closed enough to $(t_0,q(t_0))$, we deduce that for all $i\in
I_{\kappa \rho}(t,q) \cap I(t_0,q(t_0))$
 $$ \partial_t g_i(t,q) + \langle \nabla_q\, g_i(t,q),w_\nu \rangle \geq \frac{\nu \delta}{2}.$$
For the indices $i\notin I_{\kappa\rho}(t,q)$, we have
$$ g_i(t,q) + h\left[ \partial_t g_i(t,q) + \langle \nabla_q\, g_i(t,q), w_\nu \rangle\right] \geq \kappa\rho -h\beta(1+|w|+\nu).$$
Finally for $i\notin I(t_0,q(t_0))$,
$$ g_i(t,q) + h\left[ \partial_t g_i(t,q) + \langle \nabla_q\, g_i(t,q), w_\nu \rangle\right] \geq \sigma -h\beta(1+|w|+\nu)-\beta \left[|t-t_0|+|q-q(t_0)|\right],$$
with 
$$\sigma:=\min_{i\notin I(t_0,q(t_0))} g_i(t_0,q(t_0)) >0.$$
We conclude that for each fixed $\nu>0$, there are $\epsilon_\nu$ and $h_\nu$ such that for every $h<h_\nu$,
$(t,q)\in U$ with $|t-t_0|+|q-q(t_0)|\leq \epsilon_\nu$ and $q\in Q(t)$ we have
$$  \forall i, \qquad g_i(t,q) + h\left[ \partial_t g_i(t,q) + \langle \nabla_q\, g_i(t,q),w_\nu \rangle\right] \geq h\frac{\nu\delta}{2},$$
which by a first order expansion in time gives:
$$  \forall i, \qquad g_i(t+h,q) + h\langle \nabla_q\, g_i(t+h,q),w_\nu \rangle \geq h\frac{\nu\delta}{2} +O_{h\to 0}(h).$$
At the cost of decreasing $h_\nu$, it comes for $h<h_\nu$,
$$  \forall i, \qquad g_i(t+h,q) + h\langle \nabla_q\, g_i(t+h,q),w_\nu \rangle \geq 0,$$
and consequently, $w_\nu\in K_h(t+h,q)$ for every $h<h_\nu$,
$(t,q)\in U$ with $|t-t_0|+|q-q(t_0)|\leq \epsilon_\nu$ and $q\in Q(t)$.\\

\mb {\bf Step 2}:
Let us now fix the parameter $\nu$. \\
Thanks to the uniform Lipschitz regularity of the maps $q_h$ and their uniform convergence towards $q$, there exists $\tilde{h}_\nu\leq h_\nu$ such that for $\epsilon\leq \epsilon_\nu/(2+2K)$ and $h\leq \tilde{h}_\nu$,
$$ t_h^k,t_h^{k+1} \in [t_0-\epsilon,t_0+\epsilon] \Longrightarrow |t_h^{k+1}-t_0|+|q_h^k-q(t_0)|\leq \epsilon_\nu.$$
From this, as $q_h^{k}\in Q(t_h^k)$ from Proposition \ref{prop:feasible_config}, the previous step (with $t=t_h^k$) gives $w_\nu\in K_h(t_h^{k+1},q_h^{k})$. Therefore, $K_h(t_h^{k+1},q_h^{k})$ being convex, we have
\be{eq:eq}
\langle u_h^k+hf_h^k-u_h^{k+1},
w_\nu-u_h^{k+1} \rangle \leq 0.
\ee
We sum up these inequalities for $k$ from $n$ to $p$, integers chosen such that $t_h^n$ is the first time step in $[t_0-\epsilon,t_0-\epsilon+h]$ and
$t_h^{p+1}$ the last one in $[t_0+\epsilon-h,t_0+\epsilon]$. First, we know that 
\be{eq:1} \left|\sum_{k}^{p} h\langle f^k,
w_\nu-u_h^{k+1} \rangle\right| \leq \left(|w_\nu|+K\right)\int_{t_0-\epsilon}^{t_0+\epsilon} F(t) dt, \ee 
with
$K:=\sup_h \|u_h \|_\infty$. We also have 
\be{eq:2} \sum_{k}^{p} \langle u_h^k-u_h^{k+1},w_\nu \rangle = \langle
u_h(t_h^n)-u_h(t_h^{p+1}), w_\nu \rangle.\ee 
We deal with the remainder as follows: we write
$$ \sum_{k}^{p} \langle u_h^k-u_h^{k+1}, -u_h^{k+1} \rangle = \sum_{k}^{p} \langle u_h^k-u_h^{k+1}, u_h^{k} \rangle - |u_h^n|^2+|u_h^{p+1}|^2,$$
which gives
\begin{align}
 \sum_{k}^{p} \langle u_h^k-u_h^{k+1}, -u_h^{k+1} \rangle & = \frac{1}{2} \sum_{k}^p |u_h^k-u_h^{k+1}|^2 + \frac{1}{2}\left[ -|u_h(t_h^n)|^2+|u_h(t_h^{p+1})|^2\right] \nonumber \\
 & = \frac{1}{2}\textrm{Var}_2(u_h)^2_{[t_h^n,t_h^p]}+ \frac{1}{2}\left[ -|u_h(t_h^n)|^2+|u_h(t_h^p)|^2\right], \label{eq:3}
\end{align}
where we wrote $\textrm{Var}_2$ for the $L^2$-variation of a function. Using (\ref{eq:eq}), (\ref{eq:1}), (\ref{eq:2}) and
(\ref{eq:3}), we finally get~:
$$ \frac{1}{2}\textrm{Var}_2(u_h)_{[t_h^n,t_h^p]}^2+ \frac{1}{2}\left[ -|u_h(t_h^n)|^2+|u_h(t_h^{p+1})|^2\right] + \langle u_h(t_h^n)-u_h(t_h^{p}), w_\nu \rangle \lesssim \int_{t_0-\epsilon}^{t_0+\epsilon} F(t) dt.$$
Let us now choose a sequence of $\epsilon_m$ going to zero, such that $u_h$ pointwisely converges to $u$ at the instants $t_0-\epsilon_m$ and $t_0+\epsilon_m$ (which is possible as $u_h$ converges almost everywhere towards $u$). For each $\epsilon_m$ and $h$ small enough, we have shown that the last inequality holds. Then, passing to the limit for $h\to 0$ we get
\begin{align*} \frac{1}{2}\textrm{Var}_2(u)_{[t_0-\epsilon_m,t_0+\epsilon_m]}^2+ \frac{1}{2}\left[ -|u(t_0-\epsilon_m)|^2+|u(t_0+\epsilon_m)|^2\right] & \\
 & \hspace{-4cm}+\langle u(t_0-\epsilon_m)-u(t_0+\epsilon_m), w_\nu \rangle\lesssim \int_{t_0-\epsilon_m}^{t_0+\epsilon_m} F(t) dt,
\end{align*}
which gives for $\epsilon_m\to 0$
$$ \frac{1}{2}\textrm{Var}_2(u)_{[t_0^-,t_0^+]}^2+ \frac{1}{2}\left[ -|u^-(t_0)|^2+|u^+(t_0)|^2\right] + \langle u^-(t_0)-u^+(t_0), w_\nu \rangle\leq 0.$$
Finally we obtain
$$ \frac{1}{2}\left| u^{+}(t_0) - u^-(t_0)\right|^2 + \frac{1}{2}\left[ -|u^-(t_0)|^2+|u^+(t_0)|^2\right] + \langle u^-(t_0)-u^+(t_0), w_\nu \rangle\leq 0.$$
By expanding the square quantities, this can be written as follows
\be{eq:amontrer_nu} \langle u^-(t_0)-u^+(t_0), w_\nu-u^+(t_0) \rangle \leq 0.\ee

\mb
To conclude the proof, it now suffices to remember that $w_\nu=w+\nu v$ and, since for each $\nu>0$, the previous reasoning holds, we obtain (\ref{amontrer}) by letting $\nu$ go to $0$ in~(\ref{eq:amontrer_nu}). \findem

\section{Application to the modelling of inelastic collisions} \label{sec:model}

\mb
{\bf The continuous model}

\begin{figure}[hbtp]
\psfragscanon
\psfrag{qi}[l]{$q_i$}
\psfrag{qj}[l]{ $q_j$}
\psfrag{eij}[l]{$e_{ij}(q)$}
\psfrag{Dij}[l]{$D_{ij}(q)$}
\begin{center}
\resizebox{0.25\textwidth}{!}{\includegraphics{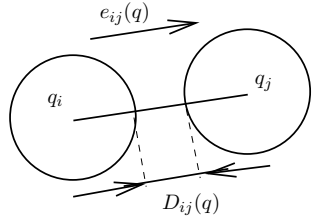}}
\end{center}
\caption{Particles $i$ and $j$ : notations.}
\label{schema_dim2_ij}
\end{figure}

\medskip
\noindent We consider a mechanical system of $N$ spherical rigid particles in three-dimensions. We denote by $q_i\in \R^3$ the position of the center of particle $i$, by $r_i$ its radius, by $m_i$ its mass and by $f_i\in \R^3$ the external force exerted on it.
Let $q\in\R^{3N}$ be defined by $q:=(\ldots,q_i,\ldots)$ and
$f\in\R^{3N}$ by $f:=(\ldots,f_i,\ldots)$. We denote by $D_{ij}(q)$ the signed distance between
particles $i$ and $j$: $$D_{ij}(q):=|q_i-q_j|-(r_i+r_j),$$
and we set $e_{ij}(q)=(q_j-q_i)/|q_j-q_i|$
(see Fig.~\ref{schema_dim2_ij}). 

\noindent The problem we are interested in is to describe the path of the configuration $q$ submitted to the force-field $f$ and undergoing inelastic collisions. This inelastic collision law can be modelled by imposing non-overlapping contraints on the particles (see the work of J.J Moreau \cite{Moreau3} introducing this concept). Therefore, we write that the positions of the particles have to belong to a set of admissible configurations $Q_0$ avoiding overlappings:
$$ q\in Q_0:= \bigcap_{i,j} \left\{q, \ D_{ij}(q)\geq 0 \right\}.$$

\mb We define $M$ as the
mass matrix of dimension $3N\times 3N$,
$M=diag(\ldots,m_i,m_i,m_i,\ldots)$. 
Then, we denote by $G_{ij}\in\R^{3N}$ the gradient of
distance $D_{ij}$ with respect to the positions of the particles:
$$
\begin{array}{cccccc}
G_{ij}(q)=&(\ldots ,0,&-e_{ij}&, 0,\ldots,0,&e_{ij}&, 0,\ldots,0)^t.\\
 & &i& &j&\\
\end{array}
$$

\mb The set $\C_q$ is the set of admissible velocities:
\be{def:Cq2}  \C_q:=\left\lbrace u,\ \langle G_{ij}(q),u\rangle\geq 0,\ \hbox{ if }\  D_{ij}(q)=0 \right\rbrace.\ee
To finish with notations, we denote by $\lambda=(\ldots,\lambda_{ij},\ldots)\in\R^{N(N-1)/2}$ the
vector made of the Lagrange multipliers associated to the $N(N-1)/2$
constraints ``$D_{ij}(q)\geq 0$''.

\mb Let $I=]0,T[$ be the time interval. The multi-particle model  we are interested in may be formally phrased as follows:
\be{Pb:cont3}
\left\{
\begin{array}{l}
\dsp q\in W^{1,\infty}(I,\R^{3N})\virg \dot q\in BV(I,\R^{3N})\virg\lambda\in ({\cal
  M}_+(I))^{N(N-1)/2}\hbox{,}\vsp\\
\dsp \forall t\in I,\quad \dot q(t^+) = \PPP_{\C_{q(t)}}\dot q(t^-)\vsp \\
\dsp M \ddot q =  f + \sum_{i<j}\lambda_{ij} G_{ij}(q) \\
\dsp \hbox{supp}(\lambda_{ij})\subset \{t\virg D_{ij}(q(t))=0\}\hbox{ for all } i,j\vsp\\
\dsp D_{ij}(q(t))\geq 0\hbox{ for all } i,j\vsp\\
\dsp q(0)=q^0 \hbox{ such that } D_{ij}(q^0)>0\hbox{ for all } i,j\virg
\dot q (0) =u^0
\end{array}
\right.
\ee

\mb The main equation 
\be{eq:main} M \ddot q - f = \sum_{i<j}\lambda_{ij} G_{ij}(q)\in -\NN(Q_0,q)\ee
 expresses the fact that overlapping is prevented by a repulsive force (the impulsion) acting on each sphere along the normal vector at the contact point. When there is no contact, $\NN(Q_0,q)$ is reduced to $\{0\}$, so that (\ref{eq:main}) reads as $M \ddot q = f$, which is the Fundamental Principle of Dynamics  applied on each sphere. Equation $\dot q(t^+) = \PPP_{\C_{q(t)}}\dot q(t^-)$ provides the inelastic collision model. It can be extended to an elastic collision model with a restitution coefficient $e$ by writing
$$ \dot q(t^+) = \PPP_{\C_{q(t)}}\dot q(t^-) - e\PPP_{\NN(Q_0,q(t))}\dot q(t^-). $$

\mb We assume for simplicity that each mass $m_i$ is equal to 1. Then Problem (\ref{Pb:cont3}) fits into the previously studied framework. 
\begin{rem}
The case of different masses can be taken into account by using the adapted scalar product $(u,v)_M=\left<Mu,v\right>$, as done in \cite{Maury}. It turns back to replace the projection step in the numerical algorithm by
$$
u^{n+1}=\PPP_{ K_h(t_h^{n+1},q_h^{n}) } \left( u^n+hM^{-1}f^n\right),
$$ 
where $\PPP$ here denotes the projection relatively to this new norm. \\
$M$ being a diagonal matrix with non-negative diagonal coefficients, it is easy to show that the following results still hold true in that case. 
\end{rem}
We emphasize that Assumption (\ref{Ui}) is satisfied as soon as 
$$ \min_{i} r_i >0,$$
and then Assumptions (\ref{gradg}) 
and (\ref{hessg}) hold true.  \\
In order to apply our previous results, it remains to check Assumption (\ref{inegtrianginverserho}). As explained in \cite{Juliette}, that corresponds to estimate the Kuhn-Tucker multipliers. Such an estimate is given in the following lemma. 
\begin{lem}
\label{lem:mult_bornes} There exists $a>0$ (depending on $N$ and on the radii $r_i$) such that
for all $q\in\R^{3N}$, $F\in\R^{3N}$ and Lagrange multipliers $(\mu_{ij})\in\R^{N(N-1)/2}$ satisfying
$$
\sum \mu_{ij} G_{ij}(q)=F \hbox{ with } \mu_{ij}\geq 0 \hbox{ and } \mu_{ij}=0 \hbox{ when } D_{ij}(q)>0,
$$
then
$$
\mu_{ij}\leq a |F|.
$$
\end{lem}
Concerning the proof of this lemma, we refer the reader to Proposition 4.7 of \cite{Juliette} (for a geometric proof) and to Proposition 2.18 of \cite{Juliette-M} (for a more ``physical'' proof). These proofs are written in a $2$-dimensional framework but they can be easily extended in our $3$-dimensional case. Actually, Lemma~\ref{lem:mult_bornes} is equivalent to Assumption (\ref{inegtrianginverserho}) with $\rho=0$. However, it can be extended and still holds for $\rho$ small enough (for example $\rho < \inf_i r_i$), see Remark 4.11 of \cite{Juliette}. 
Consequently, Assumption (\ref{inegtrianginverserho}) is satisfied for some small enough $\rho>0$.

\mb According to our main Theorem (Theorem \ref{thm}), it follows that Problem (\ref{Pb:cont3}) has solutions and that the associated numerical scheme converges (up to a subsequence). We can allow the radii to depend on time as soon as $r_i$ is uniformly twice-differentiable in time and
$$ \inf_{t\in[0,T]} \inf_{i} r_i(t) >0.$$
These theoretical results permit to legitimate the implementation of this numerical scheme. This was performed by the second author by creating SCoPI Software \cite{SCoPI}. We refer the reader to \cite{Maury} for some good properties of stability and robustness for the algorithm and efficiency for large time steps.

\begin{rem} We refer the reader to \cite{Aline}, where the second author extends this model in order to consider gluey particles. In this case, she add an extra parameter (depending on $q$) for describing the corresponding admissible set. This new operation does not keep the necessary regularity of the admissible set. She has already obtained a result of convergence for the associated numerical scheme in the single-constraint case. We plan in a forthcoming work to extend this proof with the ideas presented here in order to deal with the multi-constraint case. 
\end{rem}

\end{document}